\documentclass[letterpaper, 10 pt, conference]{ieeeconf} 
\IEEEoverridecommandlockouts 
\usepackage{graphicx}
\usepackage{subfigure}
\usepackage{amsmath}
\usepackage{amsfonts}
\usepackage[psamsfonts]{amssymb}
\usepackage{comment}
\usepackage{times}
\usepackage{changepage}
\usepackage{cases}
\usepackage{mathtools}
\usepackage{color}

\newtheorem{lemma}{Lemma}
\newtheorem{corollary}{Corollary}

\newtheorem{theorem}{Theorem}

\newcommand*{\bfrac}[2]{\genfrac{}{}{0pt}{}{#1}{#2}}

\IEEEtriggeratref{4}

\newcommand{\R}{\mathbb{R}}

\newcommand{\bmu}{\boldsymbol{\mu}}
\newcommand{\by}{\boldsymbol{y}}
\newcommand{\bz}{\boldsymbol{z}}

\newcommand{\sm}{\hspace{-1mm}-\hspace{-1mm}}
\newcommand{\spp}{\hspace{-1mm}+\hspace{-1mm}}

\title{\LARGE \bf
Cloud-Based Optimization: 
A Quasi-Decentralized Approach to Multi-Agent Coordination
}

\author{M.T. Hale and M. Egerstedt$^\dag$\thanks{$^\dag$The authors are with the School of Electrical and Computer Engineering, Georgia Institute of
Technology, Atlanta, GA 30332, USA. Email: \texttt{\{matthale, magnus\}@gatech.edu}. Research supported in
part by  the NSF under Grant CNS-1239225.}
}

\begin{document}

\maketitle
\thispagestyle{empty}
\pagestyle{empty}

\begin{abstract}
New architectures and algorithms are needed to reflect the mixture of local and global information that is available as multi-agent systems connect over the cloud. 
We present a novel architecture for multi-agent coordination where the cloud is assumed to be able to gather information from all agents, perform centralized 
computations, and disseminate the results in an intermittent manner.
This architecture is used to solve a multi-agent 
optimization problem in which each agent has a local objective
function unknown to the other agents and in which 
the agents are collectively subject to global inequality constraints.
Leveraging the cloud, a
dual problem is formulated and solved by finding a
saddle point of the associated Lagrangian. 
\end{abstract}

\section{Introduction}
Distributed optimization and algorithms have received significant attention during the last decade, e.g., \cite{bertsekas89,simon91,hendrickson00,cortes14,nedic13,johansson13,droge14},
due to the emergence of a number of application domains in which individual decision makers have to collectively arrive at a decision in a distributed manner. 
Examples of these applications include communication networks \cite{kelly98,chiang07}, sensor networks \cite{khan09,trigoni13,bullo14}, multi-robot systems \cite{cassandras13,rus13}, 
and smart power grids \cite{caron10}.

Distributed algorithms are needed mainly because the scale of large distributed systems is such that no central, global decision 
maker can collect all relevant information, perform all required computations, and then disseminate the results back to individual 
nodes in the network in a timely fashion. However, one can envision a scenario in which such globally obtained information can be used in conjunction 
with local computations performed across the network. This could, for example, be the case when a cloud computer is available 
to collect information, as was envisioned in \cite{goldberg13}. 
The question then becomes that of designing the appropriate architecture and algorithms that can leverage this mix of 
prompt decentralized computations with intermittent centralized computations.

One approach to multi-agent optimization that will prove useful towards achieving this hybrid architecture is based on 
primal-dual methods to find saddle points
of a problem's Lagrangian \cite{samar07,feijer10}. 
In fact, the study of saddle point dynamics in optimization can be traced back to
earlier results from Uzawa in \cite{arrow58}, which will provide
the starting point for the work in this paper. The primary difference between this paper and the
established literature is the cloud-based architecture
used to solve the problem; indeed the architecture is this paper's
main contribution. 
The architecture we introduce uses a cloud computer 
in order to receive information from each agent, perform global
computations, and transmit this information 
to other agents. 
We will see that this
division of labor results in globally asymptotic convergence
to an $\epsilon$-ball about a Lagrangian's saddle point. 

The goal of this paper is to serve
as a first attempt at understanding how centralized, cloud-based 
information might be injected in an intermittent but useful
manner into a network of agents where such information would 
otherwise be absent. 
In order to highlight how the cloud might
prove useful to such a system, we choose
to consider an extreme
case where no inter-agent communication occurs at all, in contrast
to existing distributed multi-agent optimization techniques, e.g., \cite{noz1, noz2, noz3, noz4}. 
Under this architecture the cloud handles all communications, and computations
are divided between the cloud and the agents in the network. 

The rest of the paper is organized as follows: Section II gives
a detailed problem statement and describes the cloud architecture,
and then Section III provides the convergence analysis for the given problem.
Next, Section IV provides
numerical results to demonstrate the viability of this approach, and
finally Section V concludes the paper.

\section{Problem Statement and Architecture}
\subsection{Architecture Motivation}
We now explain the interplay between the cloud architecture and the problem under consideration here. 
A detailed explanation is given below, with a summary and example following at the end of this section.
Consider a collection of $N$ agents indexed by $i \in A$, ${A = \{1, \ldots, N\}}$, 
where each agent is associated with a scalar state $x_i \in \R$ and where there
is no communication at all between the agents. Let the task agent $i$ is trying to solve be encoded in a strictly convex 
objective functions in $C^2$, $f_i : \R \to \R$. Each agent is assumed to have no knowledge of other agents' objective
function and each agent's only goal is to minimize its own objective function. 

To that end, agent $i$ is assumed to have immediate access to its own state, which seemingly makes this problem very simple. However, 
what prevents agent $i$ from simply computing $\frac{df_i}{dx_i}$ and setting this equal to zero -- a completely decentralized operation 
as $f_i$ only depends on $x_i$ -- is that the agents need to coordinate their actions through a globally defined constraint, 
that can, for example, represent finite resources that must be shared across the team. 
In this paper
it is assumed that agent $i$ cannot measure the state of any other agents and, as mentioned above,
that there is no communication between agents. 
Instead, this information must be obtained in some other manner, which is where the cloud will enter into the picture.

The team-wide coordination is encoded through the global constraint
\begin{equation}
g(x) = \left(\begin{array}{c}
                            g_1(x) \\
                            g_2(x) \\
                            \vdots \\
                            g_m(x) \end{array}\right) \leq 0,
\end{equation}
where $x = (x_1, \ldots, x_N)^T$ is a state vector containing the states of all agents in the network.
It is further assumed that each $g_j(x) \in C^2$ is convex. The cloud architecture discussed
here applies to any problem in which the user has selected functions $f_i$ and $g_j$ that meet the
above criteria and the forthcoming analysis fully characterizes all such problem formulations.  

Let 
\begin{equation}
F(x) = \sum_{i=1}^{N} f_i(x_i).
\end{equation}
Then $F$ is strictly convex and the problem under consideration becomes that of minimizing $F$ subject to $g$. 
The Kuhn-Tucker Theorem on concave programming (e.g., \cite{uzawa58}) states that
the optimum of this constrained problem is a saddle point of the 
Lagrangian 
\begin{equation} \label{eq:bigl}
L(x, \mu) = F(x) + \mu^{T}g(x),
\end{equation}
 where the Kuhn-Tucker (KT) multipliers $\mu_j$ satisfy $\mu_j \geq 0$ for all 
${j \in \{1, \ldots, m\}}$. 
We assume that the minimizer
of $L$ with respect to $x$, denoted $\hat{x}$, is a regular point of $g$ so 
that there is a unique saddle point, $(\hat{x}, \hat{\mu})$, of $L$ \cite{chachuat07}.
Using that $L$ is convex in $x$ and concave in $\mu$, the saddle point $(\hat{x}, \hat{\mu})$
can be shown to satisfy the inequalities
\begin{equation}
L(\hat{x}, \mu) \leq L(\hat{x}, \hat{\mu}) \leq L(x, \hat{\mu})
\end{equation}
for all admissible $x$ and $\mu$. 

Using Uzawa's algorithm \cite{arrow58}, the problem
of finding $(\hat{x}, \hat{\mu})$ can be solved from the initial point $(x(0), \mu(0))$ using the difference equations
\begin{equation} \label{eq:sysx}
x(k) = x(k-1) - \rho\frac{\partial L}{\partial x}(x(k-1), \mu(k-1))
\end{equation}
\begin{equation} \label{eq:sysmu}
\mu(k) = \max\left\{0, \,\, \mu(k-1) + \rho\frac{\partial L}{\partial \mu}(x(k-1), \mu(k-1)) \right\}
\end{equation}
where $\rho > 0$ is a constant, and where the maximum defining $\mu$ is taken component-wise so
that each component of $\mu$ is projected onto the non-negative orthant of $\R^m$, denoted by
$\R^m_+$. In the context of Uzawa's algorithm, the $i^{th}$ element of the 
state vector $x$ is updated according to 
\begin{equation} \label{eq:uzawabad}
x_i(k) = x_i(k-1) - \rho\frac{\partial L}{\partial x_i}(x(k-1), \mu(k-1)).
\end{equation}

Under the envisioned organization of the agents and the lack of inter-agent communication, 
Uzawa's algorithm cannot be directly applied. To see this, observe that if agent $i$ is to compute its own state update
using Equation
\eqref{eq:uzawabad}, 
a fundamental problem is encountered: computing $\frac{\partial L}{\partial x_i}$ 
will require knowledge of states of (possibly all) other agents and agent $i$ cannot directly 
access this information. Furthermore, determining $\mu$ at each timestep using \eqref{eq:sysmu} will 
also require the full state vector $x$, which no single agent has direct access to.

To account for the need of each agent for global information in applying
Equation \eqref{eq:uzawabad} and to compute $\mu$
using aggregated global information, the cloud computer is used. 
The cloud computer is taken to be capable of large batch computations
and receives periodic
transmissions from each agent containing
each agent's own state. 
The cloud computer uses the agents' states to compute the
next value of $\mu$ using Equation \eqref{eq:sysmu} and then transmits the states it 
received and the newly computed $\mu$ vector  to each agent. Each agent then uses 
the information from the cloud to update its own state in the vein of \eqref{eq:uzawabad}. 

\subsection{Formal Architecture Description} \label{subsec:formal}
We first describe the actions taken to initialize the system and then 
explain its operation. 
Let the agents each be programmed with their objective functions onboard and let them either
be programmed with an initial state or else be able to sense it (e.g., if it corresponds
to some physical quantity). The agents are assumed to be identifiable according
to their indices in $A$ so that the cloud knows the source of each transmission
it receives. Each agent stores and manipulates a state vector
onboard and we denote the state vector stored onboard agent $i$
by $x^i$; agent $i$'s copy of its own state is denoted $x^i_i$ and
when we are referring to a specific point in time, say timestep $k$,
we denote agent $i$'s copy of its own state at this time by $x^i_i(k)$. 
The vector of KT multipliers stored onboard agent $i$ at time $k$
is denoted $\mu^i(k)$, though we emphasize that agent $i$ does not compute
any KT vectors but instead relies on the cloud for these computations. 

Before the optimization process begins, let agent $i$ send its initial state, $x^i_i(0)$,
to the cloud and let the cloud store these states in the vector $x^c(0) \in \R^N$, with
the superscript 'c' denoting ``cloud'' and the timestep $0$ reflecting that
this is the initial state. In this notation, the cloud's copy of agent $i$'s
state at time $k$ is denoted $x^c_i(k)$. Similarly, we denote
the KT vector stored in the cloud at time $k$ by $\mu^c(k)$. 
Let the cloud be programmed by the user
with the constraint functions, $g(x^c)$. Upon receiving the
each agent's state, the cloud symbolically
computes $\frac{\partial g}{\partial x^c_i}$ and sends this function to 
agent $i$ along with some initial KT multiplier vector, $\mu(0)$,
a stepsize $\rho > 0$,
and the vector $\by^i \in \R^{N-1}$ defined as
\begin{equation} \label{eq:yidef}
\by^i = \left(\begin{array}{c}
                    x^c_1 \\
                     \vdots  \\
                     x^c_{i-1} \\
                     x^c_{i+1} \\
                     \vdots \\
                    x^c_N \end{array}\right).
\end{equation}
This vector contains states stored by the cloud in the vector $x^c$ and contains information
originally from time $0$ (though in Equation \eqref{eq:yidef} 
explicit timesteps are intentionally omitted). The subscripts in \eqref{eq:yidef} denote that agent $i$ does not
receive its own old state value from the cloud, which is logical since agent $i$ always knows its own
most recent state. In $\by^i$, then,  the cloud sends to agent $i$ the most recent state information
it has about each \emph{other} agent. 
In this notation, agent $j$'s state in $\by^i$ is denoted $\by^i_j$. In the forthcoming analysis, $\by^i$ 
always refers to the most recent state information that agent $i$ has 
received from the cloud and it will not be written as an explicit function of any time step. Similarly, 
the notation $\bmu^i$ refers to the most recent KT vector sent to agent $i$ and will be written
without an explicit timestep. 
We use the notation $\bz^i$ to denote the most recent transmission to agent $i$ containing both $\by^i$ and $\bmu^i$. 

After receiving $\bz^i$ for this first time, all $N$ agents and the cloud have the same information onboard,
and each agent begins
the optimization process. 
At timestep $0$, each agent takes one gradient step to update its own state
according to Equation \eqref{eq:uzawabad}. Simultaneously, and also at
timestep $0$, the cloud takes one gradient step to update the KT
multipliers in the cloud according to Equation \eqref{eq:sysmu}. 
Then at timestep $1$, agent $i$ sends its state, $x^i_i(1)$, to the cloud.
These transmissions are received at timestep $2$. In timestep $2$,
the cloud sends $\by^i$ and $\mu^c(1)$ to agent $i$.
These vectors are received in timestep $3$ 
at which point the cloud updates $\mu^c$ as before and each
agent takes a step to update its own state as before, thus repeating this cycle
of communication and computation. 
It is important to note that 
communications cycles do not overlap and that the agents do not send their states to the cloud at every 
timestep, but instead do so every $3^{rd}$ timestep. In addition, we emphasize
that each agent's objective function is assumed to be private throughout this process.

Due to the communications structure of the system, it is often
the case that $x^i_i(k) \neq x^j_i(k)$, namely that agents $i$ and $j$ will 
have different values for agent $i$'s state beacuse agent $j$ must wait to received agent $i$'s
state from the cloud. Due these differences, Equation \eqref{eq:sysx} is modified to reflect that each agent stores and 
manipulates a local copy of the problem. The global system therefore contains $N$
copies of the system in Equations \eqref{eq:sysx} and \eqref{eq:sysmu} and the
state vector of agent $i$ at time $k$, $x^i(k)$, is assumed to be different from that of agent $j$
at time $k$, $x^j(k)$, when $i \neq j$. 

Using the fact that agent $i$ will only update its state in timesteps just after it
receives an update from the cloud,
Equation \eqref{eq:sysx} is modified so that onboard agent $i$ it is

\begin{subnumcases} {\hspace{-5mm} \label{eq:modx} x^i(k)  = \hspace{-1mm} }
\bar{\by}^i - \rho \nabla^iL(\bar{\by}^i, \bmu^i) & \hspace{-7.5mm} $\bz^i$ received at time $k \sm 1$ \\
x^i(k-1) & else,
\end{subnumcases}
where we define 
\begin{equation} \label{eq:modgrad}
                           \nabla^i L(\bar{\by}^i, \bmu^i) = \left(\hspace{-3mm} \begin{array}{c} 0 \\ 
                                                    \vdots \\ 
             \frac{df_i}{dx_i}(x^i_i(k-1)) + (\bmu^i)^T\frac{\partial g}{\partial x_i}(\bar{\by}^i) \\ 
                                                    \vdots \\ 0 
\end{array}\hspace{-3mm} \right).
\end{equation}
Here, 
$\bar{\by}^i$ is defined as a vector onboard agent $i$ which contains $\by^i$ and the most recent
state of agent $i$ inserted in the appropriate place. In essence, $\bar{\by}^i$ is the most up-to-date
information about all of the agents that agent $i$ has access to and contains the correct value
of each other agent's state when it is received. 
Note that 
$\nabla^i L$ is simply $\frac{\partial L}{\partial x}$ with all entries except the $i^{th}$ set to $0$. 
This is because agent $i$
 does not itself compute any updates for the other agents' states which it stores onboard, but instead waits for the
 cloud to provide such updates. 

Under the architecture of this problem, only the cloud computes values of $\mu$ and there is therefore only a single update equation needed for $\mu$. 
Bearing in mind that updates to $\mu$ are only made in timesteps immediately after those in which the cloud receives each agent's state,
Equation $\eqref{eq:sysmu}$ is modified to take the form
\begin{subnumcases} {\label{eq:modmu}\hspace{-5mm} \mu^c(k)  =  \hspace{-2mm}}
\hspace{-2mm}\left[\mu^c(k\sm1) \spp \rho\frac{\partial L}{\partial \mu}(x^c(k\sm1), \mu^c(k\sm1))\right]_{+} & \hspace{-8mm} $\bfrac{\textnormal{update}}{\textnormal{at } k-1}$ \label{eq:modmureal} \\
\mu^c(k-1) & \hspace{-10mm} else, \label{eq:modmufake}
\end{subnumcases} 
where $[\cdot]_+$ denotes the projection onto $\R^m_+$ and the update referred to in Equation \eqref{eq:modmureal} is an update
of each state's value sent to the cloud. 

We note that Equation \eqref{eq:modmu} is not indexed on a per-agent basis since only the cloud computes values of $\mu$. However, we will
continue to use the notation $\mu^i(k)$ to denote the $\mu$ vector stored on agent $i$ at time $k$ (which may be different from the $\mu$ vector stored in the cloud
at time $k$). It is important to note that the argument of $\mu^i(k)$ is intended to reflect the time at which agent $i$ has $\mu^i$ onboard
and \emph{does not} imply that $\mu^i$ was computed at time $k$ or that agent $i$ computed it. 
In this notation $\bmu^i$ represents the $\mu$ vector most recently sent from the cloud to agent $i$, while 
$\mu^i$ represents the $\mu$ vector stored on agent $i$.

With this model in mind, instead of considering the system defined in Equations \eqref{eq:sysx} and \eqref{eq:sysmu}, we consider $N$ copies
of the system defined by \eqref{eq:modx} and \eqref{eq:modmu}. Using the notation that $\mu^i(k)$ represents
the vector $\mu$ as stored on agent $i$ at time $k$, we can write the full update equations onboard agent $i$ as
\begin{subnumcases} { \label{eq:agentix} \hspace{-8mm}  x^i(k) \hspace{-1mm} = \hspace{-2mm} }
\bar{\by}^i - \rho \nabla^iL(\bar{\by}^i, \bmu^i) & \hspace{-7mm} $\bz^i$ received at time $k \sm 1$ \label{eq:agentix1}\\
x^i(k-1) & else, \label{eq:agentix2}
\end{subnumcases}
\begin{subnumcases}{ \label{eq:agentimu} \mu^i(k) = }
                                                  \bmu^i & \hspace{-15mm} $\bz^i \textnormal{ was received at time } k - 1$ \\
                                                  \mu^i(k-1) & else,
\end{subnumcases}
where all changes in $\mu$ will result from the cloud using Equation \eqref{eq:modmu}.


To illustrate the communications cycle described above, Table 1 contains a sample schedule
for a single cycle. Each timestep is listed on the
left and the corresponding actions taken at that timestep are listed on the right. 


\begin{table}[h!]
\centering
\vspace{2mm}
\begin{tabular}{|c| p{2.5in} |}
\hline
Timestep & Actions \\ \hline \hline
$k$   &  Each agent receives a transmission from the cloud and then
         takes $1$ step in its own copy of the problem to update (only) its own state using Equation \eqref{eq:agentix1}. 
         At the same time, the cloud computes updated $\mu$ values using \eqref{eq:modmureal}.  \\ \hline  
$k+1$ &  Each agent sends it state to the cloud. Equation \eqref{eq:agentix2} is used by the agents and Equation \eqref{eq:modmufake} is used
         by the cloud so that no further computations are carried out during this timestep. \\ \hline
$k+2$ &  The cloud receives the agents' transmissions from time $k+1$ and stores them in $x^c$. It then sends $\by^i$ to agent $i$, along
         with $\mu^c(k+1)$, the most recently computed vector of KT multipliers (computed in timestep $k+1$). As in timestep $k+1$, Equations \eqref{eq:agentix2}
         and \eqref{eq:modmufake} are used so that no further computations take place across the network.\\ \hline
$k+3$ &  This step is identical to step $k$. Agent $i$ receives $\bz^i$ and then takes $1$ step in its own copy of the problem to update (only) its own state using Equation \eqref{eq:agentix1}. 
         At the same time, the cloud computes updated $\mu$ values using \eqref{eq:modmureal}.  \\ \hline  
\end{tabular} \par
\bigskip
Table 1: A sample schedule for one communications cycle used by the agents and cloud to exchange information.
\end{table}

\section{Convergence Analysis}
\subsection{Ultimate Boundedness of Solutions}
In this section we will examine the evolution of the sequence
\begin{equation} \label{eq:zdef}
z^i(k) = \left(\begin{array}{c} x^i(k) \\ \mu^i(k) \end{array}\right)
\end{equation}
for an arbitrary $i \in A$ in order to show that each agent's local 
copy of the problem converges to an $\epsilon$-ball
about the point $\hat{z} = (\hat{x}, \hat{\mu})$. 


Specifically, the goal here
is two-fold: to prove that the state of each agent's optimization problem
enters
a ball of radius $\epsilon$ about the saddle point $\hat{z}$
in finite time
and to show that it does not leave that ball thereafter. 
Our approach will differ from that of
\cite{arrow58} because we use the notion of ultimate
boundedness, published after Uzawa, to simplify certain components of proof.
 We restate the definition of ultimate
boundedness here for general discrete-time systems of the form
\begin{equation} \label{eq:diff}
w(k) = f(w(k-1)).
\end{equation}

\begin{lemma} \label{lem:ub}
Let $G \subseteq \R^N$ and let $V(w)$ be a Lyapunov candidate function
for the system in Equation \eqref{eq:diff}
defined on $G$ such that for all $w \in G$
\begin{equation} 
\Delta V(w) = V(f(w)) - V(w) \leq a
\end{equation}
for some $a \geq 0$. Let $\bar{G}$ denote the closure of $G$ and let $S$
be the set
\begin{equation}
S = \{w \in \bar{G} : \Delta V(w) \geq 0\}.
\end{equation}
Let $b = \sup \{V(w) : w \in S\}$ and define the set $A$ by
\begin{equation}
A = \{z \in \bar{G} : V(w) \leq a + b\}.
\end{equation}
Then any solution $\{w(k)\}$ to Equation \eqref{eq:diff} which remains in $G$ for all time
and enters $A$ at some point is contained in $A$ for all time thereafter. 
\end{lemma}
\emph{Proof:}
See \cite{hurt67}, Theorem 5. \hfill $\blacksquare$

We also state a corollary to this result which will be used below. 

\begin{corollary} \label{cor:conv}
Let the conditions of Lemma \ref{lem:ub} hold. Suppose that
\begin{equation}
\sup \{-\Delta V(w) : w  \in \bar{G} \backslash A\} > 0
\end{equation}
and that $G$ is of the form
\begin{equation}
G = \{w : V(w) \leq r\}.
\end{equation}
Then every solution $\{w(k)\}$ of Equation \eqref{eq:diff} which starts in $G$ remains in $G$
for all time and enters $A$ in a finite number of steps.
\end{corollary}
\emph{Proof:} See \cite{hurt67}, Corollaries 3 and 4. \hfill $\blacksquare$

Below, we will combine Lemma \ref{lem:ub} and Corollary \ref{cor:conv} to show that
each agent's state trajectory enters a ball of radius $\epsilon$ about $\hat{z}$,
denoted $B_{\epsilon}(\hat{z})$, and stays within that ball. 
Before proving the main convergence result, we prove the following
lemma which establishes a positive upper bound on the stepsizes that can be used. 
We will proceed in the vein of \cite{arrow58} and consider 
the Lyapunov function
\begin{equation} \label{eq:bigvdef}
V(x, \mu) = \|x - \hat{x}\|^2 + \|\mu - \hat{\mu}\|^2.
\end{equation}
\begin{lemma} \label{lem:rho}
Let $L$ denote the Lagrangian in Equation \eqref{eq:bigl} and set $R = \max\{\epsilon, V(x(0), \mu(0))\}$. Define the constants $\gamma_1$ and $\gamma_2$ by
\begin{equation}
\gamma_1 = \min_{(x, \mu)}\left\{\sqrt{\frac{\epsilon/2}{\left\|\frac{\partial L}{\partial x}\right\|^2 + \left\|\frac{\partial L}{\partial \mu}\right\|^2}} \Bigg| V(x, \mu) \leq \frac{\epsilon}{2} \right\}
\end{equation}
and
\begin{equation}
\gamma_2 = \min_{(x, \mu)}\left\{\frac{-(\hat{x} - x)^T\frac{\partial L}{\partial x} + (\hat{\mu} - \mu)^T\frac{\partial L}{\partial \mu}}{\left\|\frac{\partial L}{\partial x}\right\|^2 + \left\|\frac{\partial L}{\partial \mu}\right\|^2} \Bigg| \frac{\epsilon}{2} \leq V(x, \mu) \leq R\right\}
\end{equation}
where $\frac{\partial L}{\partial x}$ and $\frac{\partial L}{\partial \mu}$ above are (implicitly) functions of any $x$ and $\mu$ satisfying the conditions on $V$ pertaining to each set. 
Then setting
\begin{align} \label{eq:case1rho}
\rho_{max} = \min\{\gamma_1, \gamma_2\}
\end{align}
provides $\rho_{max} > 0$. 
\end{lemma}
\emph{Proof:} 
It suffices to show that $\gamma_1$ and $\gamma_2$ are both positive. The
denominator of $\gamma_1$ is always positive and tends to zero as $(x, \mu) \to (\hat{x}, \hat{\mu})$
so that the minimum defining $\gamma_1$ does not go to zero at $(\hat{x}, \hat{\mu})$. The numerator of $\gamma_1$ is
positive by insepction and $\gamma_1$ itself is therefore the square root of a positive real number.

For $\gamma_2$, we note that $L$ is convex in $x$ and concave in $\mu$. The
term $-(\hat{x} - x)^T\frac{\partial L}{\partial x}$
is the negation of the directional derivative of $L(\cdot, \mu)$ with respect
to $x$ pointing toward its minimizer, and the term $(\hat{\mu} - \mu)^T\frac{\partial L}{\partial \mu}$
is the directional derivative of $L(x, \cdot)$ with respect to $\mu$ pointing toward its maximizer. Both terms
are therefore non-negative and because the definition of $\gamma_2$ precludes $(x, \mu) = (\hat{x}, \hat{\mu})$,
the sum of these two terms is strictly positive. The denominator in the definition of $\gamma_2$ is positive
as well so that $\gamma_2$ itself is.  \hfill $\blacksquare$

The above Lemmata and Corollaries are stated in terms of the Lagrangian defined in Equation \eqref{eq:bigl}. 
While each agent in the network stores and manipulates its own state vector and thus has its own (unique) Lagrangian,
after each gradient descent step is taken and all states and KT multipliers are shared across the network,
every agent ends up with the same information before taking its next step. In addition, every agent and the cloud
use the same stepsize, $\rho$. Then despite the distribution of information and
computation throughout the network, the effective outcome of each cycle of communication and computation
as described in Section \ref{subsec:formal} is one step in each of Equations \eqref{eq:sysx} and \eqref{eq:sysmu} performed
simultaneously. 

Therefore, the analysis of the algorithm can be carried out for Equations \eqref{eq:sysx} and \eqref{eq:sysmu}, and for simplicity we 
choose to use Equations 
\eqref{eq:sysx} and \eqref{eq:sysmu} in the forthcoming analysis with the understanding that it applies equally well to all agents. 
Due to the centrality of the convergence of Uzawa's algorithm to this paper and
in order to make use of results published after the algorithm's original publication,
we now present the main result on the ultimate boundedness of solutions to the problem
at hand. 

\begin{theorem}
Let every agent use a strictly convex objective function $f_i : \R \to \R$, $f_i \in C^2$ and
let the global constraints, $g: \R^N \to \R^m$, be convex with $g_j \in C^2$ for each $j$. 
Then for any stepsize $\rho$ such that  
$0 < \rho \leq \rho_{max}$ used by all agents and the cloud,
each agent's local copy of the problem enters an $\epsilon$-ball about
$\hat{z}$ in a finite number of steps and stays within that ball
for all time thereafter. 
\end{theorem}

\emph{Proof:}
In addition to using $V$ as defined in Equation \eqref{eq:bigvdef}, 
we equivalently use that $z(k) = (x(k), \mu(k))^T$ to write
\begin{equation} \label{eq:vdef}
V(z(k)) = \|z(k) -\hat{z}\|^2.
\end{equation}
When it is convenient, we will also use the more concise notation $V(k) = V(z(k))$. 

We further define
\begin{align}
\Delta V(k) &= \|z(k+1) - \hat{z}\|^2 - \|z(k) - \hat{z}\|^2 \\
              &= \left(\|\mu(k+1) - \hat{\mu}\|^2 - \|\mu(k) - \hat{\mu}\|^2\right) \\
              &\qquad + \left(\|x(k+1) - \hat{x}\|^2 - \|x(k) - \hat{x}\|^2\right).
\end{align}
As in Lemma \ref{lem:rho}, we define
\begin{equation}
R = \max\left\{\epsilon, V\big(x(0), \mu(0)\big)\right\}.
\end{equation}

Let gradient steps be taken at timestep $k$ so that Equation \eqref{eq:agentix1} is used by
all agents to update their states and Equation \eqref{eq:modmureal} is used by the cloud to update $\mu$. 
From Equation \eqref{eq:sysx} we see that
\begin{multline} \label{eq:pfx1}
\|x(k+1)\|^2 = \|x(k)\|^2 \\ - 2\rho x(k)^T\frac{\partial L}{\partial x}(x(k), \mu(k)) + \rho^2\left\|\frac{\partial L}{\partial x}(x(k), \mu(k))\right\|^2,
\end{multline}
and multiplying both sides of Equation \eqref{eq:sysx} by $-2\hat{x}^T$ gives
\begin{equation} \label{eq:pfx2}
-2\hat{x}^Tx(k+1) = -2\hat{x}^Tx(k) - 2\hat{x}^T\frac{\partial L}{\partial x}(x(k), \mu(k)).
\end{equation}
Using Equations \eqref{eq:pfx1} and \eqref{eq:pfx2} we see that
\begin{multline} \label{eq:pfbigx}
\|x(k+1) - \hat{x}\|^2 = \|x(k) - \hat{x}\|^2 \\ + 2\rho(\hat{x} - x(k))^T\frac{\partial L}{\partial x}(x(k), \mu(k)) + \rho^2\left\|\frac{\partial L}{\partial x}(x(k), \mu(k))\right\|^2. 
\end{multline}

Carrying out the same steps for $\mu$ using Equation \eqref{eq:sysmu} gives
\begin{multline} \label{eq:pfbigmu}
\|\mu(k+1) - \hat{\mu}\|^2 = \|\mu(k) - \hat{\mu}\|^2 \\ - 2\rho(\hat{\mu} - \mu(k))^T\frac{\partial L}{\partial \mu}(x(k), \mu(k)) + \rho^2\left\|\frac{\partial L}{\partial \mu}(x(k), \mu(k))\right\|^2.
\end{multline}
Summing Equations \eqref{eq:pfbigx} and \eqref{eq:pfbigmu} gives
\begin{multline}
\|x(k+1) - \hat{x}\|^2 + \|\mu(k+1) - \hat{\mu}\|^2 \\ = \|x(k) - \hat{x}\|^2 + \|\mu(k) - \hat{\mu}\|^2 \\ - \rho\Bigg[2\bigg(-\big(\hat{x} - x(k)\big)^T\frac{\partial L}{\partial x} + \big(\hat{\mu} - \mu(k)\big)^T\frac{\partial L}{\partial \mu}\bigg) \\ - \rho\left(\left\|\frac{\partial L}{\partial x}\right\|^2 + \left\|\frac{\partial L}{\partial \mu}\right\|^2\right)\Bigg],
\end{multline}
and hence
\begin{multline} \label{eq:vdelt}
\Delta V(k) = \\ - \rho\Bigg[2\bigg(-\big(\hat{x} - x(k)\big)^T\frac{\partial L}{\partial x} + \big(\hat{\mu} - \mu(k)\big)^T\frac{\partial L}{\partial \mu}\bigg) \\ - \rho\left(\left\|\frac{\partial L}{\partial x}\right\|^2 + \left\|\frac{\partial L}{\partial \mu}\right\|^2\right)\Bigg].
\end{multline}

Suppose now that $\frac{\epsilon}{2} \leq V(k) \leq R$. Then using the fact that $\rho \leq \gamma_2$ we see that
\begin{equation} 
\Delta V(k) \leq \rho\bigg[\big(\hat{x} - x(k)\big)^T\frac{\partial L}{\partial x} - \big(\hat{\mu} - \mu(k)\big)^T\frac{\partial L}{\partial \mu}\bigg] < 0
\end{equation}
where the right-hand side is negative because $\rho$ is positive and the term inside brackets is negative.
The negativity of the term in the brackets is established by observing that it is the numerator of 
the term defining $\gamma_2$ multiplied by $-1$ and, because the numerator of the fraction defining $\gamma_2$ was shown to be positive,
we see here that this term is negative. 
In fact, the term in brackets is bounded above by some negative constant, i.e., there exists $\delta > 0$ such that
\begin{equation}
\big(\hat{x} - x(k)\big)^T\frac{\partial L}{\partial x} - \big(\hat{\mu} - \mu(k)\big)^T\frac{\partial L}{\partial \mu} \leq -\delta < 0
\end{equation}
which is seen to be true because the additive inverse of this term was shown to be bounded below by a positive constant when $\gamma_2$ was defined. 
Then for any $k$ satisfying $\frac{\epsilon}{2} \leq V(k) \leq R$, we see that $\Delta V(k) \leq -\rho\delta$ for some $\delta > 0$. 

Now suppose that $V(k) \leq \frac{\epsilon}{2}$. Then using Equation \eqref{eq:vdelt} and the fact that $\rho \leq \gamma_1$ we see that
\begin{equation}
\Delta V(k) \leq \frac{\epsilon}{2}.
\end{equation}
Here we see that $\Delta V(k) \geq 0$ only for $z(k) \in B_{\epsilon/2}(\hat{z})$ and that $\Delta V(k) \leq \frac{\epsilon}{2}$ in the set $B_{\epsilon/2}(\hat{z})$. 
Then the conditions of 
Lemma \ref{lem:ub} are satisfied with $a = b = \frac{\epsilon}{2}$ and $A = B_{\epsilon}(\hat{z})$.
In addition, for Corollary \ref{cor:conv} we see that for
any $z(k)$ satisfying $\frac{\epsilon}{2} \leq V(z(k)) \leq R$, 
there is some $\delta > 0$ such that $\sup\{-\Delta V(k)\} \geq \rho\delta > 0$. Moreover, the
set $G$ takes the form $\{z : V(z) \leq R\}$. Then
the conditions of Corollary \ref{cor:conv} are satisfied as well. 
Then $z(k)$ enters $B_{\epsilon}(\hat{z})$ in a finite
number of steps and does not ever leave thereafter. 
\hfill $\blacksquare$

To summarize, a radially unbounded, discrete-time Lyapunov function was constructed. 
The Lyapunov function was shown to satisfy the conditions needed for ultimate
boundedness and the system's trajectory was shown to come within $\epsilon$
of the Lagrangian's saddle point in finite time
and never to be more than $\epsilon$ away thereafter. 

\subsection{Extension to Private Optimization}
While above only each agent's objective function is assumed to be private,
we can extend this problem to the case where individual states are kept private. 
To do this, we modify the initialization of the system. When the cloud sends
to agent $i$ the function $\frac{\partial g}{\partial x^c_i}$, rather
than initializing agent $i$ with $\frac{\partial g}{\partial x^c_i}$ as a function of, e.g., $(x_2, x_6, x_7)$,
it can instead initialize agent $i$ with $\frac{\partial g}{\partial x^c_i}$ as a function
of $\eta^i = (\eta_1, \eta_2, \eta_3)$, where, unbeknownst to agent $i$,
$\eta_1 = x_2$, $\eta_2 = x_6$, and $\eta_3 = x_7$. By hiding
the labels of each state which will be later sent to agent $i$, these
states are kept private in the sense that agent $i$ does not know which
agent they belong to.


\section{Simulation Results}

\begin{figure}
\centering
\includegraphics[width=3.2in]{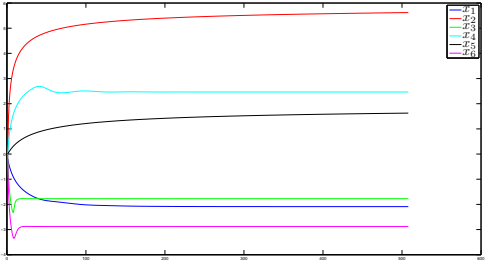}
\caption{The states onboard agent $1$ over time. Since this is a gradient-based method with a fixed stepsize, we see
larger changes in earlier iterations, followed by smaller steps taken at later iterations. 
}
\label{fig:x1states}
\end{figure}

\begin{figure}
\centering
\includegraphics[width=3.2in]{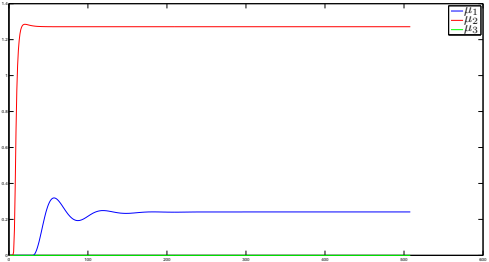}
\caption{The Kuhn-Tucker multipliers onboard agent $1$ over time. As with the states, we see
larger changes generally coming earlier in the time-evolution of the problem because they are computed using a gradient-based method
using fixed stepsizes.
}
\label{fig:x1muplot}
\end{figure}

\begin{figure}
\centering
\includegraphics[width=3.2in]{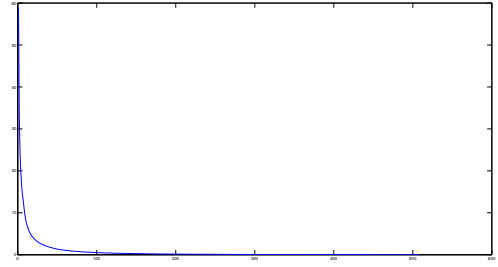}
\caption{The value of $V(x^1(k), \mu^1(k))$ over time. As was proven in the Lyapunov analysis in Section III, the Lyapunov function is non-increasing
over time.
}
\label{fig:x1vplot}
\end{figure}

A numerical implementation of the above cloud architecture was run for a particular
choice of simulation example. The problem simulated was chosen to use $N = 6$ agents,
each associated with a scalar state as above.
The objective function of each agent was chosen to be $f_i(x_i) = (x_i - \tilde{x}_i)^4$,
where
\begin{equation}
\tilde{x} = \left(\begin{array}{r} -3.0 \\ 6.0 \\ -5.0 \\ 4.0 \\ 2.0 \\ -6.0 \end{array} \right).
\end{equation}
The constraints in this problem were chosen to be
\begin{equation}
g(x) = \left(\begin{array}{c}
                    3x_1^2 + x_4^4 - 50 \\[3pt]
                    x_3^6 + x_6^4 - 100 \\[3pt]
                    9x_2 + x_5^6 - 100 \end{array}\right) \leq 0.
\end{equation}
The Lagrangian of the full problem is
\begin{equation}
L(x, \mu) = \sum_{i = 1}^6 f_i(x_i) + \mu^Tg(x)
\end{equation}
where $\mu \in \R_+^3$. 

For this example, $\gamma_1$ was found to be approximately $0.003799$ and
$\gamma_2$ was found to be approximately $0.001968$. Accordingly,
the stepsize
used was $\rho = 0.0017$. The gradient descent algorithm described above was initialized
with all agents and the cloud having all states set to $0$. All agents and the cloud
had all Kuhn-Tucker multipliers 
initialized to $0$ as well. 
Here the value $\epsilon = 0.3$ was chosen. 

For the purposes of analyzing and verifying the algorithm
presented here, the points $\hat{x}$ and $\hat{\mu}$ were computed ahead of time to be
\begin{equation}
\hat{x} = \left(\begin{array}{r}
                 -2.1278 \\
                  5.7178 \\
                 -1.7745 \\
                  2.4566 \\
                  1.6395 \\
                 -2.8798 
\end{array} \right)
\end{equation}
and
\begin{equation}
\hat{\mu} = \left(\begin{array}{c}  
      0.2462 \\
      1.2718 \\
        0  \end{array}\right).
\end{equation}

The cloud algorithm was run for $50,000$ total
iterations. It took $1,524$ iterations 
to enter a ball of radius $\epsilon$ about $\hat{z}$,
of which $508$ were spent
taking gradient descent steps and $1,016$ were spent communicating
values across the network.

The value of $x^c(50,000)$ was
\begin{equation}
x^c = \left(\begin{array}{r} 
  -2.0887 \\
   5.6219 \\
  -1.7744 \\
   2.4649 \\
   1.6271 \\
  -2.8799 
       \end{array} \right)     
\end{equation}
and the final value of $\mu^c(50,000)$ was
\begin{equation}
\mu^c = \left(\begin{array}{c} 
    0.24158 \\
    1.27176 \\
    0.00000 \\
           \end{array}\right).
\end{equation}
The final value of $V$ in the cloud was $V(x^c(50,000), \mu^c(50,000)) = 0.0110$. Based on the definition of $V$, this means that
the square of the Euclidean distance from $(x^c(50,000), \mu^c(50,000))$ to $(\hat{x}, \hat{\mu})$ is just $0.0110$.
This result confirms both that $z^c(k)$ comes within $\epsilon$ of $\hat{z}$ in finite time and that it does not go
more than $\epsilon$ away from $\hat{z}$ after that. 


To further illustrate the convergence of this problem, the histories
of the states, Kuhn-Tucker multipliers,
and value of $V$ over time onboard agent 1 for all $50,000$ timesteps are shown in Figures 1, 2, and 3, respectively. 
That $V$ is non-increasing in time was verified numerically in the MATLAB implementation and is evident in graph shown in Figure 3.

\section{Conclusion}
We presented a cloud architecture for coordinating a team
of mobile agents in a distributed optimization task. Each agent
has direct knowledge only of its own local objective function and its own
influence upon the global constraint functions but receives occasional
updates from the cloud computer containing values of each other
agent's state and updated Kuhn-Tucker multipliers. 
Using this architecture, inequality constrained multi-agent
optimization problems were proven to come within $\epsilon$
of the constrained minimum in finite time and to never
be more than $\epsilon$ away thereafter.
 Simulation results
were provided to attest to the viability of this approach.

\bibliographystyle{plain}{}
\bibliography{sources}

\end{document}